\documentclass{article}
\usepackage{amsfonts,amssymb,amsmath}

\usepackage[cp1250]{inputenc}

\usepackage{graphicx}
\usepackage{algorithm}
\usepackage{algorithmic}
\usepackage{amsmath, amsthm, amssymb}
\usepackage{amsfonts}
\usepackage{enumerate}
\usepackage{verbatim}
\usepackage{color}

\newcommand{\tluste}[1]{\mbox{\mathversion{bold}$ #1 $}}
\newcommand{\A}[0]{{\tluste{A}}}
\newcommand{\B}[0]{{\tluste{b}}}
\newcommand{\X}[0]{{\tluste{x}}}
\newcommand{\Y}[0]{{\tluste{y}}}

\newcommand{\ol}[1]{\mbox{$\overline{{#1}}$}} 
\newcommand{\ul}[1]{\mbox{$\underline{{#1}}$}}

\newcommand{\www}{\mathop{\rm w}\nolimits}
\newcommand{\intwidth}{\mathop{\rm W}\nolimits}
\newcommand{\intvolume}{\mathop{\rm V}\nolimits}

\begin{document}
\begin{center}

{\LARGE\bf
Subsquares Approach -- Simple Scheme for Solving Overdetermined Interval Linear Systems 
}\bigskip

{\large
Jaroslav Hor\'a\v cek$^{\, *}$ and Milan Hlad\'ik$^{ \, *}$
}
\vspace{1em}

$^*${ \small
 Charles University, Faculty of Mathematics and Physics, Department of Applied
Mathematics, Malostransk\'e n\'am. 25, 118 00, Prague, Czech Republic
{\tt
horacek@kam.mff.cuni.cz, hladik@kam.mff.cuni.cz
}
}
\end{center}
\vspace{0.5cm}
\noindent {\small {\bf Abstract.}
In this work we present a new simple but efficient scheme -- Subsquares approach -- for development of algorithms for enclosing the solution set of overdetermined interval linear systems. We are going to show two algorithms based on this scheme and discuss their features. We start with a simple algorithm as a motivation, then we continue with a sequential algorithm. Both algorithms can be easily parallelized. The features of both algorithms will be discussed and numerically tested.}\\


\section{Introduction}
In this paper we address the problem of solving overdetermined interval linear systems (OILS). They can occur in many situations e.g. computing eigenvectors of interval matrices \cite{hladik:eigenval} or when solving various continuous CSP problems.
There exist a lot of efficient methods for solving square interval linear systems. Solving overdetermined systems is a little bit more tricky, that is because we can not use some favourable properties of matrices like diagonal dominance, positive definiteness etc. Nevertheless there are some methods --  Rohn method \cite{rohn:enclosing}, the least squares approach \cite{arxivhoracek:oils}, linear programming \cite{arxivhoracek:oils}, Gaussian elimination \cite{hansen:ols} or the method designed by Popova \cite{popova:ols}.

In our text \cite{arxivhoracek:oils} we showed that one of the best methods is using the least squares approach. This method returns a very narrow enclosure in a small time. But there is a problem, that the least squares always return solution, even if the system has none. That is sometimes not desirable. Other methods often rely on some kind of preconditioning which leads also to an overestimation and for some systems (e.g. with really wide intervals) can not be done. It is very difficult to develop one method suitable for every type of systems. We would like to present a scheme -- Subsquares approach -- which enables us to develop methods for solving overdetermined interval linear systems. Then we will move to a simple method and sequential method, both suitable for parallel computing. Before introducing the scheme and derived methods, it would be desirable to start with some basic interval notation and definitions first.

\section{Basics of interval arithmetics}
We will work here with real closed intervals $\X = [\ul{x}, \ol{x}]$, where $\ul{x} \leq \ol{x}$. The numbers $\ul{x}, \ol{x}$ are called the lower bound and upper bound respectively.

We will use intervals as coefficients of matrices and vectors during our computations. The interval representation may be useful in many ways -- it may represent uncertainty (e.g. lack of knowledge, damage of data), verification (e.g errors in measurement), computational errors (e.g. rounding errors in floating point arithmetic) etc. Intervals and interval vectors will be denoted in boldface i.e. $\tluste{x}, \tluste{b}$. Interval matrices will be denoted by bold capitals i.e. $\A,\tluste{C}$.

%
Another notion we will use is the \emph{midpoint} of an interval $\X$, it is defined as
$x_c = (\ul{x} + \ol{x} ) / 2$. By $A_c$ we will denote the midpoint matrix of $\A$. When comparing two intervals we need the notion \emph{width} of an interval $\X$ defined as $\www(\X) = \ol{x} - \ul{x}$. If $\tluste{u}$ is an $n$-dimensional interval vector we will define "width" and "volume" of $\tluste{u}$ as
$$\intwidth(\tluste{u}) = \sum^n_{i=1} \www(\tluste{u}_i), \quad \intvolume(\tluste{u}) = \prod^n_{i=1} \www(\tluste{u}_i).$$

The vector and matrix operations will be defined using the standard interval arithmetic,  for definition of interval operations and further properties of the interval arithmetic do not hesitate to see e.g. \cite{moore:introduction}. 

We continue with definitions connected with interval linear systems. Let us have an interval linear system $\A x = \B$, where $\A$ is an $m \times n$ matrix. When $m=n$ we will call it a \emph{square} system. When $m > n$ we will call it an \emph{overdetermined} system. In the further text when talking about an overdetermined system, we will always use the notation $\A x = \B$, where $\A$ is $m \times n$ matrix. 
  
It is necessary to state what do we mean by the solution set of an interval linear system. It is the set
$$\Sigma = \{ \, x \ | \ Ax=b \ \textrm{for some} \ A \in \A, b \in \B \, \}.$$ 
 We shall point out, that this approach is different from the least squares method.
 
If a system has no solution, we call it \emph{unsolvable}. The \emph{interval hull} is an $n$-dimensional box (aligned with axes) enclosing the solution set as tightly as possible.
When we start using intervals in our computations (we have mentioned its advantage already), many problems become NP-hard \cite{kreinovich:np}.  
So is the problem of finding the hull of the solution set \cite{rohn:np}. It can be computed quite painfully using e.g. linear programming \cite{arxivhoracek:oils}. That is why we are looking for a little wider $n$-dimensional box containing the hull. The tighter the better. We call it interval \emph{enclosure}.

In this work we will provide numerical testing at various places in the text, therefore we rather mention its parameters here. The testing will be done on CPU Intel T2400, Core Duo, 1.83 GHz, with 2.50 GB memory. We used Matlab R2009a and toolbox for interval computation INTLAB v6 \cite{rump:intlab} and Versoft v10 \cite{rohn:versoft} for verified interval linear programming.

All examples will be tested on random overdetermined systems. A random system is generated in the following way. First we generate random solvable point overdetermined system. Coefficients are taken uniformly from interval $[-20, 20]$. Then we inflate the coefficients of this systems to intervals of certain width. The original point system is not necessarily a midpoint system of the new interval system. Each of its coefficients is randomly shifted towards one of the bounds of an interval in which it lies.

\section{General Subsquares method scheme}
\begin{figure}[h]
\centering 
\includegraphics[width=7cm]{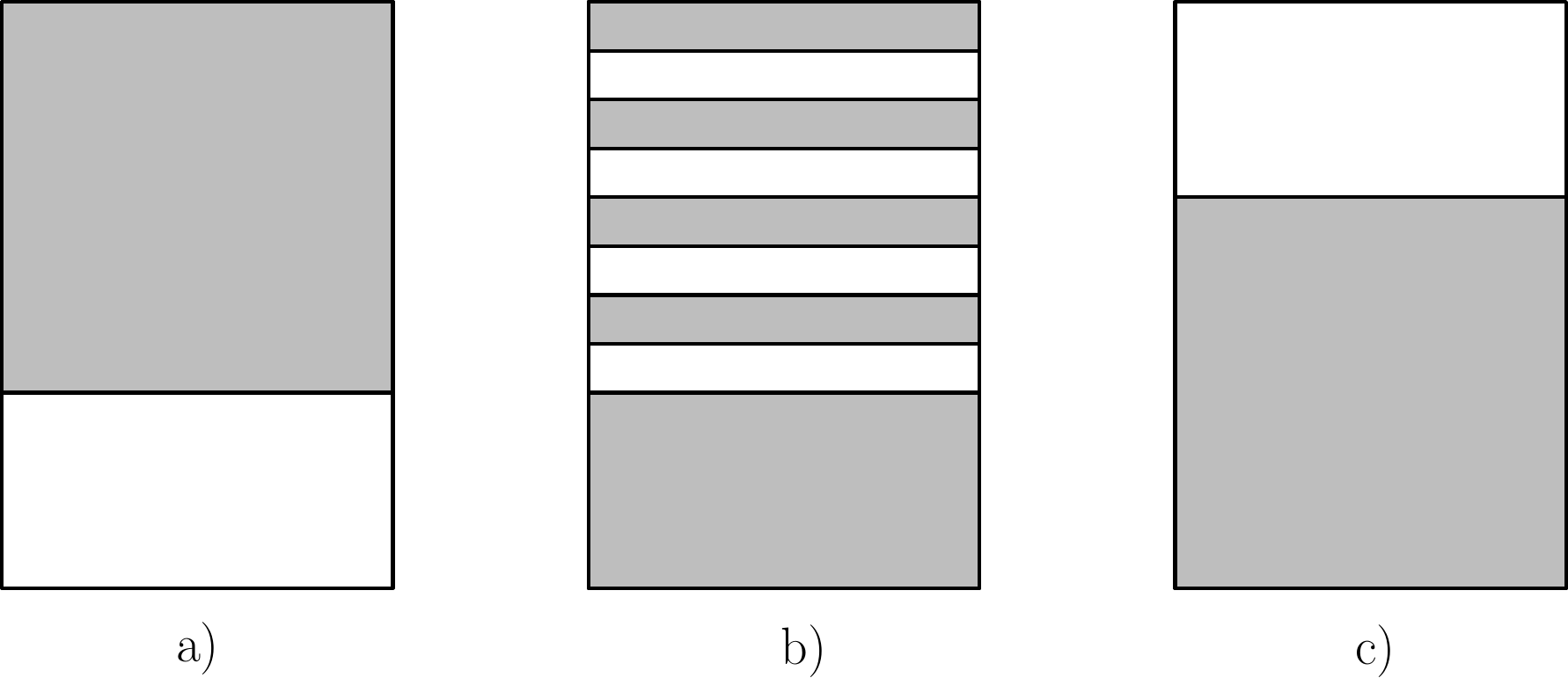}
\caption{Various square subsystems a), b), c)}
\label{fig:various_subsystems}
\end{figure}

By a square subsystem (we will also call it a \emph{subsquare}) of an overdetermined system we mean a system composed of some equations taken (without repetition) from the original overdetermined system such that together they from a square system. Some of the possibilities are shown in the Figure \ref{fig:various_subsystems}. For the sake of simplicity we will denote the square subsystem of $\A x = \B$ created by equations $i_1, i_2, \ldots, i_n$ as $\A_{\{i_1, i_2, \ldots, i_n \}} x = \B_{\{i_1, i_2, \ldots, i_n \}}$. When we use some order (e.g. dictionary order) of systems (here it does not depend which one) the $j$-th system will be denoted $\A_j x = \B_j$.

Let us suppose we can solve a square interval system efficiently and quickly. 
We can for example take one of the following method -- Jacobi method \cite{moore:introduction}, Gauss-Seidel method \cite{moore:introduction, neumaier:interval}, Krawczyk method \cite{moore:introduction, neumaier:interval} etc. These methods usually can not be applied to overdetermined systems. Nevertheless, we can use the fact that we can solve the square systems efficiently together with the fact that the solution set of an overdetermined system must lie inside the solution set of its any subsquare. This follows from the fact that by adding new equations to the square system we can only make the solution set equal or smaller.
 
\begin{figure}[h]
\centering 
\includegraphics[width=4cm]{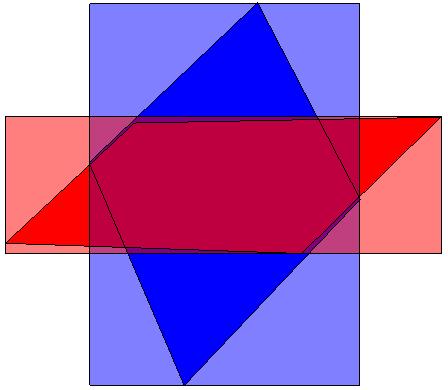}
\caption{Enclosures and hulls of subsquares}
\label{fig:enclosure_intersect}
\end{figure}

When we chose some subsquares of an overdetermined system we can simply   
provide an intersection of their solution enclosures or provide some further work. 
We get a simple scheme for solving overdetermined interval linear systems --
\emph{Subsquares Approach} --  consisting of two steps:
\begin{enumerate}
\item Select certain amount of square subsystems of $\A x = \B$
\item Solve these subsystems and get together the enclosures
\end{enumerate}
If a method for solving OILS uses this kind of approach, we call it \emph{Subsquares method}.
As a motivation for this approach let us take the randomly generated  interval system $\A x = \B$ (with rounded bounds), where
\begin{displaymath} 
\small
\A = 
\left [
\begin{array}{rrrlrrrl}
$[$&-0.8,& 0.2 &$]$ & $[$&-20.1,& -19.5&$]$\\
$[$&-15.6,& -15.2&$]$ & $[$&14.8,& 16.7&$]$\\
$[$& 18.8,& 20.1& $]$ & $[$&8.1,& 9.5& $]$\\
\end{array} \right ],
\end{displaymath} 
\begin{displaymath} 
\small
\B =  \left [
\begin{array}{rrrl}
$[$&  292.1,&  292.7&$]$\\
$[$& -361.9,& -361.1&$]$\\
$[$& 28.4,&  30.3&$]$\\
\end{array} \right ]. 
\end{displaymath} 
In the Figure \ref{fig:enclosure_intersect} we can see the solution set and the hull of $\A_{\{1,2\}} x = \B_{\{1,2\}}$ (red color) and the same for $\A_{ \{2,3 \}} x = \B_{\{2,3\}}$ (blue color). It can be seen that if we provide intersection of the two hulls (or enclosures) the resulting enclosure might get remarkably tighter. 

\subsection{The simple algorithm}

If we solve square subsystems separately and then intersect resulting enclosures, we get a raw simple Algorithm \ref{alg:simple} for solving OILS.

\begin{algorithm}
\small
\caption{Subsquares method -- simple algorithm}
\label{alg:simple}
\begin{algorithmic}
\REQUIRE $\A, \B$
\ENSURE enclosure $\X$ of the solution set of $\A x = \B$
\STATE {$\X = [-\infty, \infty]^n$}
\WHILE {\NOT (terminal condition)}
\STATE{ choose randomly a subsquare of $\A x = \B$ }
\STATE {compute its enclosure $\X_{subsq}$ } 
\STATE $\X := \X  \cap \X_{subsq} $
\ENDWHILE
\end{algorithmic}
\end{algorithm}

This approach is a little bit naive, but it has its advantage. First, if we compute enclosures of all possible square subsystems, we get really close to the interval hull. The Table \ref{tab:linvssubsq} shows the average ratios of widths and volumes of enclosures $\X_{subsq}, \ \X_{ver}$ returned by simple subsquares method and \verb|verifylss| compared to the interval hull $\X_{hull}$ computed by linear programming. 
If we have an $m \times n$ system, the number of all square subsystems is equal to $ m \choose n $.
However, we can see that for $n$ small or for $n$ close to $m$ the number $m \choose n$ might not be so large. That is why solving all subsquares pays off when systems are noodle-like or nearly-squared.

\begin{table}[h]
\centering
\small
\begin{tabular}{|c|c|c|c|c|}
\hline 
system  & av$\left( \frac{\intwidth(\X_{subsq})}{\intwidth(\X_{hull})} \right )$ & av$ \left( \frac{\intvolume(\X_{subsq})}{\intvolume(\X_{hull})}\right)$ & av$\left(\frac{\intwidth(\X_{ver})}{\intwidth(\X_{hull})}\right)$& 
av$\left(\frac{\intvolume(\X_{ver})}{\intvolume(\X_{hull})}\right)$\\[1mm]
\hline 
$5 \times 3$  & 1.0014 & 1.0043 & 1.1759 & 1.6502\\
$ 9 \times 5$ & 1.0028 & 1.0140 & 1.1906 & 2.3831 \\
$ 13 \times 7$ & 1.0044 & 1.0316 & 1.2034 & 3.6733\\
$ 15 \times 9$ & 1.0061 & 1.0565 & 1.1720& 4.2902\\
$ 25 \times 21$ & 1.0227 & 1.6060 & 1.0833 & 5.4266\\
$ 30 \times 29$  & 1.0524 & 5.8330& 1.0987 & 51.0466\\
\hline
\end{tabular}
\caption[LP vs Subsquares]{Simple subsq. method solving all subsquares -- enclosures comparison}
\label{tab:linvssubsq}
\end{table}

The second advantage is that Algorithm \ref{alg:simple} can, in contrast to other methods, easily decide whether a system is unsolvable -- if, in some iteration, the resulting enclosure is empty after intersection, then the overdetermined system is unsol\-vable. The table 
\ref{tab:subsquarenon} shows average number of random subsquares solved until the unsolvability is discovered (empty intersection occurs). Each column represents systems of different coefficient radii. We can see that for systems with relatively small intervals unsolvability is revealed almost immediately. 

\begin{table}[h]
\centering
\small
\begin{tabular}{|c|ccc|}
\hline
 system  & $rad=0.01$  & $rad=0.001$  & $rad=0.0001$ \\
\hline
$ 15 \times 10$ & 2.1 & 2.0 & 2.0\\
$ 25 \times 21 $ & 2.2 & 2.0& 2.0\\
$ 35 \times 23$ & 2.2 & 2.0 & 2.0 \\
$ 50 \times 35 $ & 2.4 & 2.0 & 2.0\\
$ 73 \times 55 $ & 2.9 & 2.1 & 2.0\\
$ 100 \times 87 $ & 7.1 & 2.1 & 2.0 \\
\hline
\end{tabular}
\caption[Subsquares -- insolvability detection]{Simple subsq. method -- unsolvability detection}
\label{tab:subsquarenon}
\end{table} 

For most rectangular systems it is however not convenient to compute solutions of all or many subsystems. The choice of square subsystems and the solving algorithm can be modified to be more economical and efficient.

\subsection{The sequential algorithm}
When selecting subsquares randomly, they usually overlap. We can think of overlaps as a "meeting points" of square subsystems. They  share some data (equations) there. That is why it would be advantageous to use this overlap to propagate a partial result of computation over one square subsystem into computations over other subsystems. When we are talking about immediate propagation of partially computed solution, our mind can easily come to Gauss-Seidel iterative method (GS). 
\newline

\noindent This method starts with an initial enclosure $\X^{(0)}$. 
In $k$-th iteration each entry of the current solution enclosure vector $\X^{(k-1)}$ might be narrowed using the formula
\begin{eqnarray}
\X^{(k)}_i & = & \frac{1}{\tluste{A}_{ii}} \Big[ \ \tluste{b}_i - ( \tluste{A}_{i1} \X^{(k)}_1 + \ldots +  \tluste{A}_{i(i-1)} \X^{(k)}_{i-1} + \nonumber \\ 
& & + \ \tluste{A}_{i(i+1)} \X^{(k-1)}_{i+1} 
+ \ldots + \tluste{A}_{in} \X^{(k-1)}_n ) \ \Big] \ \cap \ \X^{(k-1)}_i. \nonumber  
\label{gaussseidelvzorec}  
\end{eqnarray}
Simply said, in each iteration this algorithm expresses $\X_i$ from $i$-th equation. It uses newly computed values immediately. 

In our algorithm we will use GS iteration in a similar way for more square subsystems simultaneously. Again, we start with some initial enclosure $\X^{(0)}$. In $k$-th iteration we provide $k$-th GS iteration step for all systems. The advantage of this technique is that we express each variable according to formulas originating from more systems. We expect the narrowing rate will be much better this way. Similarly as in simple GS, if in some point of computation empty intersection occurs, whole overdetermined system has no solution.

Iterative methods usually require a preconditioning. We will use the preconditioning with $A_c^{-1}$.  
There are still two not answered problems yet -- initial enclosure and terminal condition. To find $\X^{(0)}$, we can take the resulting enclosure of some other method. Or we can solve one square subsystem.
The algorithm will terminate after $k$-th iteration if e.g. 
$$ \forall \, i \quad \left |\, \ul{\X}^{(k)}_{\,i} - \ul{\X}^{(k-1)}_{\,i} \right | < \epsilon \quad \textrm{and} \quad \left |\, \ol{\X}^{\,(k)}_{\,i} - \ol{\X}^{\,(k-1)}_{\,i} \right | < \epsilon,$$
for some small positive $\epsilon$ and $i=1,\ldots,n$.

Now we have to choose subsystems for sequential method.
Before getting on, we shall consider the following desirable properties of a new algorithm inspired with four unfavorable features of the simple algorithm:
\begin{enumerate}
\item We do not want to have too many square subsystems
\item We want to cover the whole overdetermined system by subsystems
\item The overlap of subsquares is not too low, not too high
\item We take subsquares that narrow the resulting enclosure as much as possible
\end{enumerate}

First and second property can be solved by covering the system step by step using some overlap parameter. About third property,  it proved itself to be a reasonable choice taking overlap $\approx n/3$. Property four is a difficult task to provide. We think deciding which systems to choose (in a favourable time) is still an area to be explored. Yet randomness will serve us well. 

Among many possibilities we tested, the following choice of subsystems worked well. During our algorithm we divide equations of overdetermined system into two sets -- $covered$, which contains equations that are already contained in some subsystems, and $waiting$, which contains equations that are not covered yet. We also use a parameter $overlap$.

The first subsystem will be chosen randomly, other subsystems will be composed of $overlap$ rows from $covered$ and $n - overlap$ rows from $waiting$. The last system is composed of all remaining uncovered rows and then some already covered rows are added to form a square system. This is described as Algorithm \ref{alg:choosingsystems}. The algorithm is not necessarily optimal, it should serve as an illustration.
The procedure \emph{randsel(n, list)} selects \emph{n} random non-repeating numbers from \emph{list}. The total number of subsquares chosen by this algorithm is $$\small 1 + \left \lceil \frac{m-n}{n-overlap}\right \rceil.$$

The whole sequential algorithm is summarized as Algorithm \ref{alg:sequential}. The function \emph{GS-iteration($\tluste{C} x=\tluste{d}$, $\Y$)} applies one iteration of Gauss-Seidel method on the subsquare $\tluste{C} x=\tluste{d}$ using $\Y$ as initial enclosure .  Method \emph{has-converged()} returns true if terminal condition (e.g. the one mentioned earlier) is satisfied.


\begin{algorithm}
\small
\caption{Choosing square subsystems}
\label{alg:choosingsystems}
\begin{algorithmic}
\REQUIRE $\A, \B$, overlap
\ENSURE set of subsquares of $\A x = \B$
\STATE {systems $\leftarrow \emptyset$} \COMMENT {set of square subsystems}
\STATE {covered $\leftarrow$ $\emptyset$}  \COMMENT {numbers of covered equations by some subsystem}
\STATE {waiting $\leftarrow$ $\{1, 2, \ldots, m\}$}   \COMMENT {numbers of equations to be covered} 
\newline
\WHILE { $\textrm{waiting} \neq \emptyset $ }
\IF {$\textrm{covered} = \emptyset$}
\STATE{indices $\leftarrow$ randsel($n$, waiting) }
\newline
\ELSIF {$|\textrm{waiting}| \leq (n - \textrm{overlap})$}
\STATE {indices $\leftarrow$ waiting $\cup$ randsel($n-|\textrm{waiting}|$, waiting)}
\newline
\ELSE
\STATE {indices $\leftarrow$ randsel(overlap, covered) $\cup$ randsel(n-overlap, waiting)}
\ENDIF
\newline
\STATE {systems $\leftarrow$ systems $\cup \ \{ \A_{\textrm{indices}} x = \B_{\textrm{indices}} \}$}
\STATE {covered $\leftarrow$ covered $\cup$ indices}
\STATE {waiting $\leftarrow$ waiting $\setminus$ indices}
\ENDWHILE
\RETURN {systems}
\end{algorithmic}
\end{algorithm}
 
\begin{algorithm}
\small
\caption{Subsquares method -- sequential version}
\label{alg:sequential}
\begin{algorithmic}
\REQUIRE $\A, \B, \X^{(0)}$
\ENSURE enclosure $\X$ of the solution set $\A x = \B$
\STATE {select square subsystems $\{\A_1 x = \B_1 , \ldots, \A_{k} x = \B_{k}\}$}
\STATE {$\X \leftarrow \X^{(0)}$ }
\STATE {converged $\leftarrow$ false}
\WHILE { \NOT converged } 
\FOR{ $i=1$ \TO $k$} 
\STATE { $\X \leftarrow$  GS-iteration$(\A_i x = \B_i, \X)$ } 
\ENDFOR
\STATE {converged $\leftarrow$ has-converged()}
\ENDWHILE
\RETURN $\X$

\end{algorithmic}
\end{algorithm}

As we showed in \cite{arxivhoracek:oils}, \verb|verifylss| (using the least squares) from INTLAB is one of the best and quite general method for overdetermined systems that are 
solvable. That is why we wanted to compare our method with this method. 
During our testing we realized \verb|verifylss| works fine with small intervals, however it is not too efficient when the intervals become relatively large. We used enclosures returned by \verb|verifylss| as inputs for subsquares sequential method and tested if our algorithm was able to narrow them. 
The Table \ref{tab:subqvsver} shows the results. Column $rad$ shows radii of intervals of testing systems, we chose the same radii for all coefficients of a system. We tested on 100 random systems. For each system we chose 100 random subsquares sets and applied the sequential algorithm on them. The fourth column shows average ratios of enclosure widths of $\X_{subsq}$ and $\X_{ver}$. Each random selection of subsquares set usually produces a different ratio of enclosure widths. For each system we chose one of the 100 subsquares sets that produces the best ratio. The fifth column shows the average value of the best ratios found for each of 100 random systems. Columns $t_{ver}$ and $t_{subsq}$ show computation times (in seconds) of \verb|verifylss| and sequential method respectively.

\begin{table}[h]
\centering
\small
\begin{tabular}{|c|c|c||c|c||c|c|}
\hline
system    & overlap & $rad $ & av$\left ( \frac{\intwidth(\X_{subsq})}{\intwidth(\X_{ver})} \right)$  & av$\left ( \textrm{best} \frac{\intwidth(\X_{subsq})}{\intwidth(\X_{ver})} \right)$  & $t_{ver}$ & $t_{subsq}$\\
\hline
$ 15 \times 10$ & 3 & 0.1& 0.99& 0.94& 0.006&  0.06\\
$ 15 \times 10$ & 3 & 0.25& 0.97& 0.86& 0.007& 0.07 \\
$ 15 \times 10$ & 3 & 0.35& 0.93 & 0.79& 0.008& 0.09 \\
$ 15 \times 10$ & 3 & 0.5& 0.87& 0.66& 0.01&0.12  \\
\hline
$ 25 \times 13 $& 5& 0.1& 0.99 & 0.98& 0.006& 0.09\\
$ 25 \times 13 $& 5& 0.25& 0.99& 0.94& 0.007& 0.12\\
$ 25 \times 13 $& 5& 0.35 & 0.98& 0.92& 0.008& 0.14\\
$ 25 \times 13 $& 5& 0.5& 0.94& 0.79& 0.012&0.20 \\
\hline
$ 37 \times 20$ & 7 &0.1 & 0.99& 0.98& 0.008& 0.11\\
$ 37 \times 20$ & 7 &0.25 & 0.99& 0.95& 0.011& 0.19\\
$ 37 \times 20$ & 7 & 0.35& 0.97& 0.90& 0.015& 0.29\\
$ 37 \times 20$ & 7 &0.5 & 0.87& 0.38& 0.016& 0.72\\
\hline
$ 50 \times 35$ & 11 & 0.1& 0.99& 0.98& 0.014& 0.16\\
$ 50 \times 35$ & 11 &0.25 & 0.97& 0.84& 0.023& 0.51\\
\hline
\end{tabular}
\caption[Subsquares vs Verifylss]{Subsquare method shaving the verifylss enclosure}
\label{tab:subqvsver}
\end{table} 
The larger interval radii are, the more intensively the sequential method sharpens \verb|verifylss| enclosures.
When an interval system has its midpoint system $(A_c x = b_c)$, or a system relatively close to it, solvable, we are able to compute the solution enclosure much tighter. This happens because of the preconditioning with $A_c^{-1}$ matrix. Selected results can be seen in the Table \ref{tab:subqvsvercenter}.   

\begin{table}[h]
\centering
\small
\begin{tabular}{|c|c|c||c|c|c|}
\hline
system    & overlap & $rad $ & av. rat  & best av. rat  \\ 
\hline
$ 15 \times 10$ & 3 & 0.1& 0.98& 0.89\\
$ 15 \times 10$ & 3 & 0.25& 0.85& 0.59 \\
$ 15 \times 10$ & 3 & 0.35& 0.76 &0.40 \\
\hline
$ 25 \times 13 $& 5& 0.1& 0.99 & 0.96\\
$ 25 \times 13 $& 5& 0.25& 0.93& 0.74\\
$ 25 \times 13 $& 5& 0.35 & 0.76& 0.33\\
\hline
$ 37 \times 20$ & 7 &0.1 & 0.99& 0.97\\
$ 37 \times 20$ & 7 &0.25 & 0.89& 0.34\\
\hline
$ 50 \times 35$ & 11 & 0.1& 0.98& 0.82\\
\hline
\end{tabular}
\caption[Subsquares vs Verifylss]{Subsquares method shaving the verifylss enclosure (solvable midpoint system)}
\label{tab:subqvsvercenter}
\end{table}

\subsection{Parallel algorithm}
The naive algorithm can be easily parallelized. All the square subsystems can be solved in parallel and then all enclosures are intersected. We have only one barrier at the end of computation. 

If we take a look at computation times in the Table \ref{tab:subqvsver}, we realize \verb|verifylss| is much faster. However, the time pay for gaining much precise enclosure using sequential method is not too high.
Moreover, even the sequential algorithm can be parallelized. Propagation of newly computed enclosure can be guaranteed by sharing the enclosure vector $\X$ among processors as a global variable. If we use Jacobi formula instead of Gauss-Seidel formula for one iteration, the computation becomes of kind SIMD - single instruction multiple data. Therefore it could be used even on GPUs -- one pipeline narrows one variable from the interval enclosure vector, a bunch of pipelines computes over one subsquare.   
Nevertheless, shared vector $\X$ might create a bottleneck. We believe this  could by prevented by the following behaviour of each pipeline. When reading, each pipeline does not lock corresponding shared variable. After each iteration it overwrites a shared variable only if it has improved the value currently stored there. This is inspired with an observation, that in one iteration not many computations over different subsquares improve the same variable. However, we assume that there exist much more efficient ways how to make parallel subsquares methods more efficient and memory collision avoiding.  

\section{Conclusion}
In this paper we introduced a simple but efficient scheme  -- subsquares approach -- for enclosing the solution set of overdetermined interval linear systems. The first method derived from this scheme was a little bit naive, but for noddle-like or nearly-square systems it was able to find almost the interval hull. The second method was a little bit more sophisticated but still quite simple. It worked well on interval systems which coefficients were composed of wide intervals. This method was able to significantly sharpen enclosures produced by \verb|verifylss|. Both methods were able to detect unsolvability of OILS. Moreover they could be easily parallelized. In the second method we chose the square subsystems randomly, that is why sharpening produced by this method had variable results. There is a question whether for each OILS there exists a deterministically chosen set of subsquares which gives the best possible enclosure, so we can avoid the rando\-mization. 
\newline

\noindent{\bf Acknowledgement.}
\newline
{ \small
Our research was supported by the grant GA\v{C}R P402/13/10660S.
Jaroslav Hor\'a\v cek was partially supported by the Grant Agency of the Charles University (GAUK) grant no. 712912 and by GAUK no. 
SVV-2013-267313. }

{
\small
\bibliographystyle{abbrv}
\bibliography{literatura}
}

\end{document}